\documentclass [12pt] {amsart}
\usepackage {amssymb, amsmath, amscd}

\theoremstyle {plain}
\newtheorem {theorem} {Theorem}

\newtheorem {lemma} {Lemma}

\theoremstyle {definition}
\newtheorem {definition} {Definition}

\theoremstyle {remark}
\newtheorem {remark} {Remark}

\newtheorem* {acknowledgments} {Acknowledgments}

\numberwithin {equation} {section}

\begin {document}
\title {Some endomorphisms of $II_1$ factors}
\author {Hsiang-Ping Huang}
\address {Department of Mathematics, University of Utah, Salt Lake City, UT 84112}
\email {hphuang@math.utah.edu}
\thanks {Research partially supported by National Center for Theoretical
         Sciences, Mathematics Division, Taiwan}

\pagestyle{plain}

\begin {abstract}
For any finite dimensional $C^*$-algebra $A$
, we
give an endomorphism $\Phi$ of
the hyperfinite $II_1$ factor $R$ of finite
Jones index such that:
\[
  \forall \ k \ \in \ \mathbb {N}, \
  \Phi^k (R)' \cap R= \otimes^k A.
\]
The Jones index $[R: \Phi (R)]= ({\rm rank\ } (A))^2$, here ${\rm  rank\ } (A)$
is the dimension of the maximal abelian subalgebra of $A$.
\end {abstract}
\maketitle
\tableofcontents
\section {Introduction} \label {S: intro}

In \cite {rP88} R.Powers initiated the study of unital
*-endomorphisms of the hyperfinite $II_1$ factor $R$.
He defined a shift $\Phi$ on $R$ to be a unital *-endomorphism of
$R$ such that $\cap_{k=1}^{\infty} \Phi^k (R)= \mathbb {C}$ and
defined the index of $\Phi$ as the Jones index $[R: \Phi (R)]$. He
constructed a family of shifts, called binary shifts, on $R$. The
index of a binary shift is 2. In \cite {mC87} M.Choda generalized
the construction and obtained a $\mathbb {N}$-parameterized family
of shifts, called n-unitary shifts; the Jones index of a n-unitary
shift is equal to n. In \cite {gP87} G.Price constructed the first
example of nonbinary shifts on $R$ of index 2. In \cite {dB88}
D.Bures and H.-S.Yin obtained an intrinsic characterization of such
shifts, called group shifts: shifts on the twisted group von Neumann
algebra.

M.Enomoto, M.Nagisa, Y.Watatani, and H.Yoshida \cite {mE91}
calculated the relative commutant algebras $\Phi^k (R)' \cap R$ of
binary shifts. Recently G.Price \cite {gP99} obtained the complete
classification of rational binary shifts. The entropy of a rational
binary shift is $\frac {1} {2}  {\log 2}$. On the other extreme of
the spectrum, \cite {hN95} H.Narnhofer, W.Thirring, and E.St{\o}rmer
utilized the so-called ''power-irreducible'' binary shifts and
obtained surprising counterexamples for which the tensor product
formula for entropy fails. In \cite {vG99} V.Golodets and
E.St{\o}rmer showed the relation of the center sequence of the
binary shift and the entropy. If the center sequence grows faster
than $0 (n)$, then the mean entropy can take any value in
$(\frac{1}{2}\log 2,\log 2]$.

In this paper we couple the notion of n-unitary shifts with the
shift on $\otimes^{\infty}_{i=1} A$, for any finite
dimensional $C^*$-algebra $A$. We then construct a family of
shifts on $R$ of integer index
and try to calculate the relative commutant algebras
and the entropy.
One of the most distinguished examples, $\Phi$ on $R$, gives:
\[
  \Phi^k (R)' \cap R= \otimes^k A.
\]
The Jones index $[R: \Phi (R)] = ({\rm  rank\ } (A))^2$, here ${\rm  rank\ } (A)$
is the dimension of the maximal abelian subalgebra of $A$.

\section {Preliminaries} \label {S: prelim}
Let $M$ be a $II_1$ factor with the canonical trace $\tau$. Denote
the set of unital *-endomorphisms of $M$ by $End (M, \tau)$. Then
$\Phi \in End (M, \tau)$ preserves the trace and $\Phi$ is
injective. $\Phi (M)$ is a subfactor of $M$. If there exists a
$\sigma \in Aut (M)$ with $\Phi_1 \cdot \sigma= \sigma \cdot \Phi_2$
for $\Phi_i \in End (M, \tau)$ $(i= 1, 2)$ then $\Phi_1$ and
$\Phi_2$ are said to be conjugate. If there exists a $\sigma \in Aut
(M)$ and a unitary $u \in M$ such that $Ad u \cdot \Phi_1 \cdot
\sigma= \sigma \cdot \Phi_2$, then $\Phi_1$ and $\Phi_2$ are outer
conjugate.

The Jones index $[M: \Phi (M)]$ is an outer-conjugacy invariant. We
consider only the finite index case unless otherwise stated. In such
case, there is a distinguished outer-conjugacy invariant: the tower
of inclusions of finite dimensional $C^*$ algebras, $\{ A_k= \Phi^k
(M)' \cap M \}_{k=1}^{\infty}$.

Another well-known conjugacy invariant is the Connes-St{\o}rmer
entropy.

\begin {lemma}
$A_k = \Phi^k (M)' \cap M$ contains an subalgebra that is isomorphic
to $\otimes_{l= 1}^k A_1$, the $k$-th tensor power of $A_1$, where
$A_1 = \Phi (M)' \cap M$ as denoted.
\end {lemma}

\begin {proof}
We collect some facts here, for $1 \leq l < k$: \\
(1) $\Phi^{l} (A_1)$ is isomorphic to $A_1$,
since $\Phi^{l}$ is injective.\\
(2) $\Phi^{l} (A_1) \cap A_1 \subseteq \Phi
(M) \cap \Phi (M)'= \mathbb{C}$.\\
(3) $[A_1, \Phi^{l} (A_1)] = 0$ by the definition of $A_1$.\\
(4) $\Phi^{l} (A_1) \subseteq A_{l+1} \subseteq A_k$.\\

Put $\{s_{{i_1, j_1}}\}_{i_1, j_1 \in F}$ to be a system of matrix
units for $A_1$. Similarly we have $\{ \Phi (s_{i_2, j_2}) \}_{i_2,
j_2 \in F}$, a system of matrix units for $\Phi (A_1)$. Etc.

Consider the following linear equation:
\[
  \sum_{i_1, j_1, \cdots, i_k, j_k} a_{i_1, j_1, \cdots, i_k, j_k}
  s_{i_1, j_1} \Phi(s_{i_2, j_2}) \cdots \Phi^{k-1} (s_{i_k, j_k})= 0.
\]

If we can conclude the coefficient $a_{i_1, j_1,\cdots, i_k, j_k}$
is zero for every $i_1, j_1,$ $\cdots,$ $i_k, j_k$, then the
dimension of $\bigvee_{l=1}^{k} \Phi^{l-1} (A_1)$ is equal to the
$k$-th power of the dimension of $A_1$. In other words, the linear
independence is established.

Multiply the equation  by $s_{i_1, i_1} \Phi(s_{i_2, i_2}) \cdots
\Phi^{k-1}(s_{i_k, j_k})$ at the left hand side and  by $s_{j_1,
i_1} \Phi(s_{j_2, i_2}) \cdots \Phi^{k-1}(s_{j_k, i_k})$ at the
right hand side, we get:
\[
  a_{i_1, j_1, \cdots, i_k, j_k} s_{i_1, i_1} \Phi (s_{i_2, i_2}) \cdots
  \Phi^{k-1}(s_{i_k, i_k})= 0.
\]

Note that $s_{i_1, i_1} \Phi(s_{i_2, i_2}) \cdots \Phi^{k-1}(s_{i_k,
i_k})$ is a projection and has nonzero trace, we can conclude:
$a_{i_1, j_1, \cdots, i_k, j_k}$ is zero for every $i_1, j_1,
\cdots, i_k, j_k$.
\[
  \bigvee_{l=1}^{k} \Phi^{l-1} (A_1) \simeq \otimes_{l=1}^{k}
  \Phi^{l-1} (A_1) \simeq \otimes_{l=1}^{k} A_1.
\]
\end {proof}

\begin {remark}
The dimension of the relative commutant ${\Phi^k (M)}' \cap M$ is
known to be bounded above by the Jones index $[M: \Phi(M)]^k$. Lemma
1 provides the lower bound for the growth estimate.
\end {remark}

A good example is the canonical shift \cite {dB} on the tower of
higher relative commutants for a strongly amenable inclusion of
$II_1$ factors of finite index. The ascending union of higher
relative commutants gives the hyperfinite $II_1$ factor, and the
canonical shift can be viewed as a *-endomorphism on the hyperfinite
factor. Lemma 1 is nothing but the commutation relations in S.Popa's
$\lambda$-lattice axioms \cite {sP95}.

A natural question arises with the above observation: for any finite
dimensional C*-algebra $A$, can we find a $II_1$ factor $M$ and a
$\Phi \in End (M, \tau)$ such that: for all $k \in {\mathbb N}$,
\[
  \Phi^k (M)' \cap M \simeq \otimes_{i=1}^k A_i,
  {\text {\ where \ }} A_i= A?
\]

The answer is positive and furthermore we can choose
$M$ to be the hyperfinite $II_1$ factor $R$.
We give the construction
in the next section. The main technical tool in
the construction is \cite {rP88} R.Powers' binary
shifts. We provide here the details of n-unitary
shifts generalized by \cite {mC87} M.Choda for the convenience
of the reader.

Let $n$ be a positive integer. We treat a pair of sets $Q$ and
$S$ of integers satisfying the following condition $(*)$ for
some integer $m$:
\[
(*)
\begin {cases}
Q= (i(1), i(2), \cdots, i(m)), &0 \leq i(1) < i(2) < \cdots
< i(m), \\
S= (j(1), j(2), \cdots, j(m)), & j(l)= 1, 2, \cdots, n-1,\\
&\text {\ for \ } l= 1, 2, \cdots, m.
\end {cases}
\]
\begin {definition}
A unital $*$-endomorphism $\Psi$ of $R$ is called an n-unitary
shift of $R$ if there is a unitary $u \in R$ satisfying the following:\\
(1)$u^n= 1$;\\
(2)$R$ is generated by $\{u, \Psi (u), \Psi^2 (u), \cdots, \}$;\\
(3)$\Psi^k (u) u= u \Psi^k (u)$ or $\Psi^k (u) u= \gamma u \Psi^k
(u)$
for all $k= 1,2,\cdots$, where $\gamma = \exp (2 \pi \sqrt {-1} / n)$.\\
(4) for each $(Q, S)$ satisfying $(*)$, there are an integer $k (\geq
0)$ and a nontrivial $\lambda \in \mathbb {T} = \{ \mu \in \mathbb {C};
| \mu | = 1 \}$ such that
\[
  \Psi^k (u) u(Q, S)= \lambda u(Q, S) \Psi^k (u),
\]
where  $u(Q, S)$ is defined by
\[
  u(Q, S) = {\Psi^{i(1)} (u)}^{j(1)} {\Psi^{i(2)} (u)}^{j(2)}
  \cdots  {\Psi^{i(m)} (u)}^{j(m)}.
\]
\end {definition}

The unitary $u$ is called a generator of $\Psi$. Put $S (\Psi; u)=
\{ k; \Psi^k (u) u= \gamma u \Psi^k (u) \}$. Note that the above
condition (2) gives some rigidity on $S (\Psi; u)$.
The Jones index $[R: \Psi (R)]$ is n.

One interesting example of $S_1= S (\Psi_1; u_1)$ is $\{
1,3,6,10,15, \cdots,$ $\frac {1} {2} l(l+1),  \cdots \}$, which
corresponds to the n-stream $\{ 0 1 0 1 0 0 1 0 0 0 1 0 0 0 0 1$ $ 0
0 0 0 0 1 \cdots \}$. It is pointed out that the relative commutant
$\Psi_1^k (R)' \cap R$ is always trivial for all $k$! That is, our
question for $A= \mathbb {C}$ is answered by this example.

In \cite {dB88}, D.Bures and H.-S.Yin introduced a notion of group
shifts, constructed by realizing $R$ as the twisted group von
Neumann algebra on a discrete abelian group with a 2-cocycle. The
special case of the group $G= \oplus_{i=1}^{\infty} \mathbb
{Z}_n^{(i)}$ with a suitable 2-cocycle and the (right) 1-shift
generalizes the above result by M.Choda. Put $\{u_1, u_2, \cdots \}$
as the set of generators of $G$. We can specify on the abelian group
$G$ any 2-cocycle $\omega$ by its associated antisymmetric character
$\rho$ of $G \wedge G$,
\begin{align*}
  &\rho (g \wedge h):= \omega (g, h) {\overline
  {\omega (h, g)}}\\
  &u_g u_h= \omega (g, h) u_{gh}= \omega (g, h) {\overline {\omega
  (h, g)}} u_h u_g= \rho (g \wedge h) u_h u_g
\end{align*}

Define a character by:
\begin {align*}
  &\rho (u_{2i-1} \wedge u_{2i})= \gamma,\\
  &\rho (u_{2i-1} \wedge u_j)= 1, \text{ if } j \not= 2i,\\
  &\rho (u_{2i} \wedge u_{2i+2j}) = \gamma, \text{ if } j \in S_1,\\
  &\rho (u_{2i} \wedge u_{2i+2j}) = 1, \text{ if } j \notin S_1.
\end {align*}
The set $S_1$ and $\gamma$ is as mentioned above.

Take the (right) 2-shift $\Psi_2$ on $G$; $\Psi_2 (u_i) = u_{i+2}$.
This 2-shift is compatible with the given 2-cocycle. Proposition 1.2
\cite {dB88} tells us that $\Psi_2$ can be extended to the twisted
group von Neumann algebra, which is the hyperfinite $II_1$ factor,
with the following property:
\[
  \Psi_2^{k} (R)' \cap R \simeq {\mathbb {C}}^{kn}; \quad
  [R: \Psi_2 (R)]= n^2.
\]

In fact, $\Psi_2^{k} (R)' \cap R$ is generated by $\{ u_1, u_3,
\cdots, u_{2k-1} \}$. Observe that the von Neumann algebra generated
by $u_{2k-1}, u_{2k}$ is isomorphic to the $n \times n$ matrix
algebra. For our question, this construction realizes the case when
$A$ is abelian. The construction for the most general case is
modeled on this special one.

\section {Main Theorem} \label {S: main}
\begin {theorem}
For any finite dimensional $C^*$-algebra $A$, there exists a $\Phi
\in End (R, \tau)$ such that the relative commutant $ \Phi^k (R)'
\cap R$ is isomorphic to $\otimes_{i=1}^k A^{(i)}$. Here $R$ is the
hyperfinite $II_1$ factor with the trace $\tau$.
\end {theorem}

Since $A$ is finite dimensional, then $A$ can be decomposed as
a direct sum of finitely many matrix algebras,
\[
  A \simeq \oplus_{i=1}^j M_{a_i} (\mathbb {C}) \subseteq M_n
  (\mathbb {C}),
\]
where $n= \sum_i a_i$.

For each $i$, $M_{a_i} (\mathbb {C}) \subset A$ (not a unital
embedding) is generated by $p_i, q_i \in \mathcal {U}(\mathbb
{C}^{a_i})$ with:
\[
  p_i^{a_i}= q_i^{a_i}= 1_{M_{a_i}(\mathbb {C})}; \quad
  \gamma_i= \exp( 2 \pi \sqrt {-1} / a_i), \quad p_i q_i = \gamma_i  q_i p_i
\]
where $p_i= [ 1 \ \gamma_i \ \gamma_i^2 \cdots \gamma_i^{a_i-1} ]$
is the diagonal matrix in $M_{a_i} (\mathbb {C})$,
 and
$q_i$ is the permutation matrix in $M_{a_i} (\mathbb {C})$,
$( 1 \ 2 \ 3 \cdots a_i )$.

Denote by $\mathcal {D}$ the diagonal algebra of $A \subseteq
M_n (\mathbb {C})$. $\mathcal {D}$ is
generated by $\{ p_i \}_{i=1}^j$. A simple observation is that $M_n
(\mathbb {C})$ is generated by $\mathcal {D}$ and a unitary
$u$ with $u^n = 1$. We can write $u$ as the permutation matrix
$(1 \ 2 \ 3 \ \cdots \ n)$.
Every element $x$ in $M_n (\mathbb {C})$ is expressed as $x=
\sum_{i=0}^{n-1} x_i u^i$, where $x_i \in \mathcal {D} \subseteq A$.

Note that $Ad u$ defines a unital *-automorphism of $\mathcal {D}$,
yet not of $A$. $Ad u$ doesn't map $A$ into $A$.
As a consequence, for $x$ in $\mathcal {D}$ we have:
\[
   u^i x = {Ad u}^i (x) u^i; \quad
   u^i \mathcal {D}= \mathcal {D} u^i.
\]

Define $v \in M_n( \mathbb {C})$ to be the permutation matrix:
\[
v= (a_1 \ \ (a_1+ a_2) \ \ (a_1+ a_2+ a_3) \ \ \cdots \ \
   (a_1+ a_2+ \cdots +a_j)).
\]
Then $v^j= 1$.
Note that $Ad v$ defines a unital *-automorphism of $\mathcal {D}$.
As a consequence, for $x$ in $\mathcal {D}$ we have:
\[
   v^i x = {Ad v}^i (x) v^i; \quad
   v^i \mathcal {D}= \mathcal {D} v^i.
\]

$\mathcal {D}$ and $v$ don't generate the full matrix algebra,
$M_n (\mathbb {C})$. However $A$ and $v$ do generate $M_n (\mathbb
{C})$. Therefore we have two ways to describe $M_n (\mathbb {C})$:\\
(1) via $\mathcal {D}$ and $u$, or\\
(2) via $A$ and $v$ (or via $A$ and $r$ described below.)\\

Not surprisingly there is a relation between $u$ and $v$;
\[
  u= v_1 v_2 v_3 \cdots v_j v,\
  \text {where} \  v_i= q_i + 1 - 1_{M_{a_i}} (\mathbb {C}).
\]
Define $r:= s v s$, while
\[
  s= [0 \ 0 \ \cdots \ 0 \ 1_{a_1} \ 0 \ 0 \
  \cdots \ 0 \ 1_{a_1+ a_2} \cdots \ 0 \ 0 \cdots \ 0 \ 1_{a_1+ a_2+ \cdots+ a_j}]
  \in \mathcal {D}.
\]
Thus $A$ and $r$ generates $M_n (\mathbb {C})$.

On the other hand, define $w= \sum_{i=1}^j \gamma^{i-1} 1_{M_{a_i}
(\mathbb {C})}$, where $\gamma= \exp (2 \pi \sqrt {-1} / j)$ and
$\gamma^j = 1$. Note that $w$ is in the center of $A$.
Two simple yet important observations are that:\\
(1) $Ad w$ acts trivially on $A$, which contains $\mathcal {D}$.\\
(2) $Ad w (r)= \gamma r$.

We now construct a tower of inclusion of finite dimensional $C^*$-
algebras $M_k$ with a trace $\tau$. The ascending union $M= \cup_{k
\in \mathbb {N}} M_k$ contains infinite copies of $M_n (\mathbb
{C})$, and thus of $A$. Number them respectively by
$\mathcal {D}_1, r_1, A_1, w_1$, $\mathcal {D}_2, r_2, A_2,
w_2$,
$\mathcal {D}_3$, $r_3$, $A_3, w_3$, $\cdots$.

We endow on this algebra the following properties:
\[
  [A_l, A_m]= 0, \ \text{if} \  l \not= m;
\]
\[
  [r_l, A_m]= 0, \ \text{if} \ l \not= m;
\]
\[
  r_l r_m = \gamma r_l r_m,\  \text{if} \  |l- m| \in
  S_1= \{ 1, 3, 6, 10, 15, \cdots \},
\]
\[
  r_l r_m = r_m r_l, \ \text{otherwise}.
\]
Here $\gamma= \exp (2 \pi \sqrt {-1} / j)$ as above.

The construction is an induction process. We have handy the
embedding $A_1 \subseteq M_n (\mathbb {C})= M_1$, which is
isomorphic to the inclusion of $A \otimes 1_{M_n (\mathbb {C})}$
inside $M_n (\mathbb {C}) \otimes 1_{M_n (\mathbb {C})}$. Identify
$A_2$ in $\otimes^2 M_n (\mathbb {C})$ by $1_{M_n (\mathbb {C})}
\otimes A$.

Consider the *-automorphism on $\otimes^2 M_n (\mathbb {C})$, $Ad (w
\otimes v)$.
We have the following results:\\
(1)$Ad (w \otimes v)$ = $Ad (1_{M_n (\mathbb {C})} \otimes v)$
\ \text {when restricted on} \ $1_{M_n (\mathbb {C})} \otimes A= A_2$.\\
(2)$Ad (w \otimes v)$ acts trivially on  $A \otimes 1_{M_n (\mathbb
{C})}= A_1$.
\\
(3)$Ad (w \otimes v) (r_1 \otimes 1_{M_n (\mathbb {C})}) =
   \gamma r_1 \otimes 1_{M_n (\mathbb {C})}$.

Take $v_2= w \otimes v \in \otimes^2 M_n ({\mathbb C})$, $v_2^j= 1$.
Define $r_2: = s_2 v_2 s_2= w \otimes r$, where $s_2$ is as above a
projection in ${\mathcal D}_2 \subset A_2$.
We have the following properties:\\
(1)$<A_2, r_2> \simeq M_n (\mathbb {C})$ by the existence of a system of matrix units; \\
(2)$[A_1, A_2]=0$,  $[A_1, r_2]= 0$;\\
(3)$Ad v_2 (r_1)= \gamma r_1$, i.e., $r_1 r_2= \gamma r_2 r_1$;\\
(4) The trace $\tau$ on $M_2= <A_1, r_1, A_2, r_2> \subset \otimes^2
M_n ({\mathbb C})$ is an extension of the normalized trace $\tau=
\frac {1} {n} Tr$ on $M_n (\mathbb {C})= M_1= <A_1, r_1>$.

\begin{lemma}
In the above construction, $M_2= <A_1, r_1, A_2, r_2>$ is equal to
$\otimes^2 M_n ({\mathbb C})$ as a concrete $C^*$-algebra.
\end{lemma}

Assume we have obtained $M_k= <A_1, r_1, A_2, r_2, \cdots, A_k,
r_k>$ equal to $ \otimes^k M_n ({\mathbb C})$ with the trace $\tau$.
We identify $M_k$ as $M_k \otimes 1_{M_n (\mathbb {C})}$ by sending
$x \in M_k$ to $x \otimes 1_{M_n (\mathbb {C})}$, and similarly
$A_{k+1}$ as $1_{M_k} \otimes A$.

Consider the *-automorphism on $M_k \otimes M_n (\mathbb {C})$, $Ad
((w_1^{b_1} w_2^{b_2} \cdots w_k^{b_k}) \otimes v)$, $\forall \ l
\leq k$:
\[
    b_l = 1, \ \text {if} \  |k+1- l| \in S_1;
  b_l= 0, \ \text {otherwise}.
\]
We have:\\
(1) $v_{k+1}^j= 1$, where $v_{k+1}$ is a unitary;\\
(2)  $r_{k+1}:= s_{k+1} v_{k+1} s_{k+1}
 = $$(w_1^{b_1} w_2^{b_2} \cdots w_k^{b_k}) \otimes r$ and the algebra generated by
 $<A_{k+1}, r_{k+1}>$ is isomorphic to $M_n (\mathbb {C})$;\\
(3) $[A_{k+1}, M_k]= 0$;\\
(4) $[r_{k+1}, A_l]= 0$, for $l \leq k$;\\
(5) $Ad v_{k+1}(r_l)= \gamma r_l$, if $|k+1- l| \in S_1$;
    $Ad v_{k+1} (r_l)= r_l$, otherwise;\\
Or equivalently to (5);\\
(5)'     $r_{k+1} r_l = \gamma r_l r_{k+1}$,
     if $|k+1- l| \in S_1$; $r_{k+1} r_l = r_l r_{k+1}$, otherwise;\\
(6) There exists an extending trace $\tau$ on the finite dimensional
$C^*$-algebra $M_{k+1}= <M_k, A_{k+1}, r_{k+1}>$.

\begin{lemma}
$M_k$ is equal to $\otimes^k M_n ({\mathbb C})$.
\end{lemma}

By induction we have constructed the ascending tower of finite
dimensional $C^*$-algebras with the desired properties.

We now explore some useful properties of the finite dimensional
$C^*$-algebra, $M_k$.

\begin {lemma}
For all $k$, $M_k$ is the linear span of the words,
$x_1 \cdot x_2 \cdot x_3 \cdot \cdots \cdot x_k$, where $x_i \in M_n
(\mathbb {C})_i= <A_i, r_i>$.
\end {lemma}

\begin {proof}
It suffices to prove $x_l x_i$ is in $<A_i, r_i> \cdot
<A_l, r_l>$, where $i < l$.
$x_l \in (A_l r_l+ A_l)^n$ and $x_i \in (A_i r_i+ A_i)^n$ through
the decomposition of $M_n (\mathbb {C})_i$ and $M_n (\mathbb {C})_l$
by $D_i, u_i, D_l, u_l$. Thus it suffices to
prove
\[
   (A_l r_l+ A_l) \cdot (A_i r_i+ A_i) \subseteq
    (A_i r_i+ A_i) \cdot (A_l r_l+ A_l).
\]
\begin {align*}
&A_l A_i= A_i A_l;\\
&A_l r_l A_i= A_i A_l r_l;\\
&A_l A_i r_i= A_i r_i A_l;\\
&A_l r_l A_i r_i= A_i A_l r_l r_i= A_i A_l r_i r_l
=A_i r_i A_l r_l.
\end {align*}
\end {proof}

\begin {lemma}
In fact, $<A, r> \simeq M_n (\mathbb {C})$ is of the form:
\[
  A+ A r A+ A r^2 A+ \cdots + A r^{j-1} A.
\]
\end {lemma}

\begin {proof}
It suffices to observe that
\begin{align*}
   &r= s v s= v s, \quad [s, v]= [s, r]=0\\
   &r^*= s v^*= v^* s, \quad r r^*= r^* r= s\\
   &r A r= r s A  r= r^2 r^* A r=  r^2 s (v^* A v) s \subset r^2 A\\
   &r^*= r^{j-1}, \quad r^j= s
\end {align*}

\end {proof}

\begin {lemma}
Consider the pair $(M, \tau)$ as described above and the GNS-
construction. Identify everything mentioned above as its image. We
$M''$ is the hyperfinite $II_1$ factor.
\end {lemma}

\begin{proof}
There is one and only one tracial state on $M_k$ for all $k \in
{\mathbb N}$. Hence the tracial state on $M$ is unique. Therefore
$M''$ is the hyperfinite $II_1$ factor, $R$.
\end {proof}

Define a unital *-endomorphism, $\Phi$, on $R$ to be the (right)
one-shift: i.e., sending $A_k$ to $A_{k+1}$, and sending $r_k$ to
$r_{k+1}$. We observe that $\Phi (R)$ is a $II_1$ factor and
\[
  [R: \Phi (R)]= n^2.
\]

\begin {lemma}
The relative commutant $\Phi^k (R)' \cap R$ is exactly
$\otimes_{i=1}^k A$.
\end {lemma}

\begin {proof}
Because of our decomposition in Lemma 4 and Lemma 5, $R$ can be
written as
\[
  (\sum_{i=0}^j A_1 r_1^{i} A_1) \cdot (\sum_{i=0}^j
  A_2 r_2^{i} A_2) \cdot \cdots \cdot (\sum_{i=0}^j A_k r_k^{i} A_k)
  \cdot \Phi^k (R).
\]

Assume $x \in R \cap \Phi^k (R)'$. $x$ can be written, as in Lemma 4,
of the following form:
\[
  x= \sum_{\vec {\alpha} \in \{ 0, \ 1,  \ \cdots, j-1 \}^k} y_1^{\vec {\alpha}} r_1^{c_1} z_1^{\vec
  {\alpha}}
  y_2^{\vec {\alpha}} r_2^{c_2} z_2^{\vec {\alpha}} \cdots y_k^{\vec {\alpha}} r_1^{c_k}
  z_k^{\vec {\alpha}}
  \cdot y^{\vec {\alpha}},
\]
where $\vec {\alpha}= (c_1, c_2, \dots, c_k)$ is a multi-index and
$y^{\vec {\alpha}}$ is in $\Phi^k (R)$. Note that $\Phi^k (R)$ is
the weak closure of $\{ \Phi^k (M_i) \}_{i=1}^{\infty}$.

For every $\epsilon > 0$,  there exists an integer $i \in {\mathbb
N}$ such that
\begin{align*}
   &\forall \ {\vec {\alpha}}, \ z^{\vec {\alpha}} \in \Phi^k (M_i)
   \subset <A_{k+1}, r_{k+1}, \cdots, A_{k+i}, r_{k+i}>\\
   &\| x-\sum_{\vec {\alpha} \in \{ 0, \ 1,  \ \cdots, j-1 \}^k} y_1^{\vec {\alpha}} r_1^{c_1} z_1^{\vec
  {\alpha}}
  y_2^{\vec {\alpha}} r_2^{c_2} z_2^{\vec {\alpha}} \cdots y_k^{\vec {\alpha}} r_1^{c_k}
  z_k^{\vec {\alpha}}
  \cdot z^{\vec {\alpha}} \|_{2, \tau} < \delta\\
  &\delta= (\sqrt{\frac{j}{n}})^k \epsilon
\end{align*}


Put $L= l( l+1)/2+1$ for some integer $l > k+1$.  We have the
following properties:
\begin {align*}
&[r_L, A_1]= [r_L, A_2]= \cdots = [r_L, A_{k+i}]= 0\\
&[r_L, r_2]= [r_L, r_3]= \cdots = [r_L, r_{k+i}]= 0\\
&r_L r_1 = \gamma r_1 r_L, \quad r_L r_1^{c_1}
   = \gamma^{c_1} r_1 r_L\\
&r_L r_L^*= r_L^* r_L= s_L\\
&r_L r_1 r_L^*= \gamma r_1 s_L, \quad r_L r_1^{c_1} r_L^*= \gamma^{c_1} r_1 s_L\\
&{\rm for \ } 0 \leq m \leq {j-1}, \quad r_L^m r_1^{c_1} (r_L^*)^m= \gamma^{c_1 m} r_1 s_l\\
&[s_L, A_1]= [s_L, A_2]= \cdots = [s_L, A_{k+i}]= 0\\
&[s_L, r_1]= [s_L, r_2]= \cdots = [s_L, r_{k+i}]= 0
\end {align*}

Therefore we claim:
\begin {align*}
&\| (x-\sum_{\vec {\alpha} \in \{ 0, \ 1,  \ \cdots, j-1 \}^k}
y_1^{\vec {\alpha}} r_1^{c_1} z_1^{\vec
  {\alpha}}
  y_2^{\vec {\alpha}} r_2^{c_2} z_2^{\vec {\alpha}} \cdots y_k^{\vec {\alpha}} r_1^{c_k}
  z_k^{\vec {\alpha}}
  \cdot z^{\vec {\alpha}}) s_L \|_{2, \tau}=\\
&\| (x-\sum_{\vec {\alpha} \in \{ 0, \ 1,  \ \cdots, j-1 \}^k}
y_1^{\vec {\alpha}} r_1^{c_1} z_1^{\vec
  {\alpha}}
  y_2^{\vec {\alpha}} r_2^{c_2} z_2^{\vec {\alpha}} \cdots y_k^{\vec {\alpha}} r_1^{c_k}
  z_k^{\vec {\alpha}}
  \cdot z^{\vec {\alpha}}) \frac{1}{j} \sum_{m=0}^{j-1} r_L^m (r_L^*)^m \|_{2, \tau}=\\
&\frac{1} {j} \| \sum_{\vec {\alpha}} \sum_m (r_L^m x {r_L^*}^m-
y_1^{\vec {\alpha}} r_L^m r_1^{c_1} (r_L^*)^m z_1^{\vec
  {\alpha}}
  y_2^{\vec {\alpha}} r_2^{c_2} z_2^{\vec {\alpha}} \cdots y_k^{\vec {\alpha}} r_1^{c_k}
  z_k^{\vec {\alpha}}
  \cdot z^{\vec {\alpha}}) \|_{2, \tau}=\\
&\frac{1} {j} \| \sum_{\vec {\alpha}} \sum_m (x - y_1^{\vec
{\alpha}}  r_1^{c_1 m} z_1^{\vec
  {\alpha}}
  y_2^{\vec {\alpha}} r_2^{c_2} z_2^{\vec {\alpha}} \cdots y_k^{\vec {\alpha}} r_1^{c_k}
  z_k^{\vec {\alpha}}
  \cdot z^{\vec {\alpha}}) s_L \|_{2, \tau}=\\
&\| (x-\sum_{\vec {\alpha} \in \{ 0, \ 1, \cdots, j-1 \}^k, c_1= 0}
y_1^{\vec {\alpha}} z_1^{\vec
  {\alpha}}
  y_2^{\vec {\alpha}} r_2^{c_2} z_2^{\vec {\alpha}} \cdots y_k^{\vec {\alpha}} r_1^{c_k}
  z_k^{\vec {\alpha}}
  \cdot z^{\vec {\alpha}}) s_L \|_{2, \tau}=\\
&\sqrt{\frac{j}{n}} \| x-\sum_{\vec {\alpha} \in \{ 0, \ 1, \cdots,
j-1 \}^k, c_1= 0} y_1^{\vec {\alpha}} z_1^{\vec
  {\alpha}}
  y_2^{\vec {\alpha}} r_2^{c_2} z_2^{\vec {\alpha}} \cdots y_k^{\vec {\alpha}} r_1^{c_k}
  z_k^{\vec {\alpha}}
  \cdot z^{\vec {\alpha}}  \|_{2, \tau}\\
&{\rm since \ } \{ x, \ y_1^{\vec {\alpha}}, \ z_1^{\vec
  {\alpha}}, \ y_2^{\vec {\alpha}}, \ r_2^{c_2}, \ z_2^{\vec {\alpha}}, \cdots, y_k^{\vec{\alpha}},
  \ r_1^{c_k}, \ z_k^{\vec {\alpha}}, \ z^{\vec {\alpha}} \} \subset \{ s_L, \ r_L, \ A_L
  \}'\\
&{\rm and \ } \text{\cite{sP83}}\tau (s_L)= \frac{j}{n} .
\end {align*}

By induction,
\begin{align*}
&\| x-\sum_{\vec {\alpha} \in \{ 0, \ 1, \cdots, j-1 \}^k, \ c_1= 0}
y_1^{\vec {\alpha}} z_1^{\vec
  {\alpha}}
  y_2^{\vec {\alpha}} r_2^{c_2} z_2^{\vec {\alpha}} \cdots y_k^{\vec {\alpha}} r_1^{c_k}
  z_k^{\vec {\alpha}}
  \cdot z^{\vec {\alpha}}  \|_{2, \tau} < \sqrt{\frac{n}{j}} \delta\\
&\| x-\sum_{\vec {\alpha} \in \{ 0, \ 1, \cdots, j-1 \}^k, \ c_1=
c_2= 0} y_1^{\vec {\alpha}} z_1^{\vec
  {\alpha}}
  y_2^{\vec {\alpha}} z_2^{\vec {\alpha}} \cdots y_k^{\vec {\alpha}} r_1^{c_k}
  z_k^{\vec {\alpha}}
  \cdot z^{\vec {\alpha}}  \|_{2, \tau} < (\sqrt{\frac{n}{j}})^2 \delta\\
&\cdots \\
&\| x-\sum_{\vec {\alpha} \in \{ 0 \}^k} y_1^{\vec {\alpha}}
z_1^{\vec
  {\alpha}}
  y_2^{\vec {\alpha}} z_2^{\vec {\alpha}} \cdots y_k^{\vec {\alpha}}
  z_k^{\vec {\alpha}}
  \cdot z^{\vec {\alpha}}  \|_{2, \tau} < (\sqrt{\frac{n}{j}})^k \delta= \epsilon\\
\end{align*}

Note that the von Neumann algebra  $\{ x, \ y_1^{\vec {\alpha}}, \
z_1^{\vec {\alpha}}, \ y_2^{\vec {\alpha}}, \ z_2^{\vec {\alpha}},
\cdots y_k^{\vec {\alpha}},\ z_k^{\vec {\alpha}} \}''$ commutes with
$\Phi^k (M)$, which is a $II_1$ factor. Any element in the former
von Neumann algebra has a scalar conditional expectation onto
$\Phi^k (M)$. In short, the former von Neumann algebra and $\Phi^k
(M)$ are mutually orthogonal.

According to \cite {sP83}, we have:
\begin{align*}
  &\| x-\sum_{\vec {\alpha} \in \{ 0 \}^k} y_1^{\vec {\alpha}}
  z_1^{\vec
  {\alpha}}
  y_2^{\vec {\alpha}} z_2^{\vec {\alpha}} \cdots y_k^{\vec {\alpha}}
  z_k^{\vec {\alpha}}
  \cdot \tau (z^{\vec {\alpha}})  \|_{2, \tau} < \epsilon\\
  &\sum_{\vec {\alpha} \in \{ 0 \}^k} y_1^{\vec {\alpha}} z_1^{\vec
  {\alpha}}
  y_2^{\vec {\alpha}} z_2^{\vec {\alpha}} \cdots y_k^{\vec {\alpha}}
  z_k^{\vec {\alpha}}
  \cdot \tau (z^{\vec {\alpha}})  \in A_1 \cdot A_2 \cdots A_k= \otimes^k A
\end{align*}

\end {proof}

\section {Discussion} \label {S:disc}

The construction of $\Phi$ depends on the choice of the
anticommutation set $S_1$. In the case of n-unitary shifts,
different choices of $S (\Psi; u)$ give uncountably many
nonconjugate shifts and at least a countably infinite family of
shifts that are pairwise not outer conjugate. Not to mention in
\cite {mC90} the existence of uncountably many non-outer-conjugate
nonbinary shifts, exploiting different 2-cocycles on the group $G=
\oplus_i^{\infty} \mathbb {Z}_2^{(i)}$. Each of the above has a
counterpart in our construction.


Given a finite dimensional $C^*$-algebra $A$ with ${\rm {rank}} (A)=
n$ and an n-unitary shift with the anticommutation set $S (\Psi; u)$
\cite {mC87} satisfies the following condition:
\[
  S( \Psi; u)= \{k_i \ | \
  k_{i+1} > k_i, \text {\ and \ } k_{i+2}- k_{i+1} > k_{i+1}- k_i,
  \ \forall
  i \in \mathbb {N} \}
\]

Denote by $M (\Psi; u)$ to be the ascending union of $M_k (\Psi;
u)=$
\[
  <A_1 (\Psi; u), r_1 (\Psi; u), A_2 (\Psi; u), r_2 (\Psi; u),
  \cdots, A_k (\Psi; u), r_k (\Psi; u)>
\]
with the trace $\tau (\Psi; u)$. We have
\[
   A \simeq A_1(\Psi; u) \simeq A_2(\Psi; u) \cdots
\]

Similarly, consider the pair $(M (\Psi; u), \tau (\Psi; u))$ as
described above and the GNS construction. Identify everything
mentioned above as its image. The weak closure $M (\Psi; u)''$ is
the hyperfinite $II_1$ factor, $R$.

Define a unital *-endomorphism, $\Psi$, on $R$ to be the (right)
one-shift: i.e., sending $A_k (\Psi; u)$ to $A_{k+1} (\Psi; u)$, and
sending $r_k (\Psi; u)$ to $r_{k+1} (\Psi; u)$. We observe that
$\Psi (R)$ is a $II_1$ factor and $[R: \Psi (R)]= n^2$.
\[
  \Psi^k (R)' \cap R= \otimes^k A.
\]

We calculate the Connes-St{\o}rmer entropy \cite {eS00} in the
following paragraph.

\begin {lemma}
$H (\Psi) \geq \ln n= \ln [R: \Psi (R)]$ no matter of the choice of
the anticommutation set $S (\Psi, u)$.
\end {lemma}

\begin {proof}
\begin {align*}
&H (\Psi)= \lim_{j \to \infty} H (M_j, \Psi) \geq \lim_{j \to
\infty}
H (\otimes^j A, \Psi)\\
&= \lim_{j \to \infty} \lim_{k \to \infty} \frac {1} {k}H (\otimes^j
A, \Psi (\otimes^j A),
\cdots, \Psi^{k-1} (\otimes^j A))\\
&= \lim_{j \to \infty} \lim_{k \to \infty} \frac {1} {k} H
(\otimes^{j +k -1} A) = \lim_{j \to \infty} \lim_{k \to \infty}
\frac {1} {k}
\sum_{n^{j+k-1}} -\frac {1} {n^{j+k-1}} \ln (\frac {1} {n^{j+k-1}}) \\
&= \lim_{j \to \infty} \ln n= \ln n
\end {align*}
\end {proof}

\begin {lemma}
$H (\Psi) \leq \ln n= \ln [R: \Psi (R)]$ no matter of the choice of
the anticommutation set $S (\Psi, u)$.
\end {lemma}

\begin {proof}
\begin {align*}
&H (\Psi)= \lim_{j \to \infty} H (M_j, \Psi) = \lim_{j \to \infty}
\lim_{k \to \infty} \frac {1} {k} H( M_j, \Psi (M_j), \cdots, \Psi^{k-1} (M_j))\\
& \leq \lim_{j \to \infty} \lim_{k \to \infty} \frac {1} {k} H(
M_{j+ k- 1})= \lim_{j \to \infty} \lim_{k \to \infty} \frac {1} {k}
H(\otimes^{j+k-1} M_n ({\mathbb C}))\\
&=\lim_{j \to \infty} \lim_{k \to \infty} \frac {1} {k} (j+k -1) \ln
n= \ln n
\end {align*}

\end {proof}

Therefore, $H (\Psi) = \ln n= \ln [R: \Psi (R)]$ no matter of the
choice of the anticommutation set $S (\Psi, u)$.

\section {Digression} \label {S:digr}
In this section, we construct inclusions of non-hyperfinite
$II_1$-factors via free products with amalgamation.

\begin {theorem}
For any finite dimensional $C^*$-algebra $A$, there exists a tower
of inclusions of $II_1$-factors, $M \subset M_1 \subset M_2 \subset
M_3 \subset \cdots$, with the trace $\tau$ such that
\[
   M' \cap M_k = \otimes_{i=1}^k A.
\]
\end {theorem}
%
%

The main tool is the relative commutant theorem by S.Popa.
\begin {lemma} \cite {sP}
Let $(P_1, \tau_1), (P_2, \tau_2)$ be two finite von Neumann
algebras with a common von Neumann subalgebra $B \subset P_1$, $B
\subset P_2$, such that $P_1= Q {\overline {\otimes}} B$ where $Q$
is a nonatomic finite von Neumann algebra. If $(P, \tau)$ denotes
the amalgamated free product $(P_1, \tau_1)*_B (P_2, \tau_1)$ then
${Q_0}' \cap P = ({Q_0}' \cap Q) \overline {\otimes} B$ for any
nonatomic von Neumann subalgebra $Q_0 \subset Q$.

Assume that $L \subset P$ is a von Neumann subalgebra satisfying
the properties:\\
(1) $Q_0 \subset L$.\\
(2) $L \cap P_2$ contains an element $y \not= 0$ orthogonal to $B$,
i.e. $E_B (y)= 0$ and with $E_B (y^* y) \in \mathbb {C} 1$.\\
Then $L' \cap P= L' \cap B$. If in addition $L \cap B= \mathbb {C}$
then $L$ is a type $II_1$ factor.
\end {lemma}

As before, we can embed $A$ in the full matrix algebra $M_n (\mathbb
{C})$. Consider the tensor product of $M_n (\mathbb {C})$ and $M_2
(\mathbb {C})$. Identify $A$ as its image in the tensor product.
Take the element $y$ in $M_2 (\mathbb {C})$:
\[
  y= {\begin {pmatrix}
       0 & 1\\
       1 & 0
      \end {pmatrix}}.
\]
In $M_{2n} (\mathbb {C})$, $y \not= 0$ is an element orthogonal to
$A$, i.e., $E_{A} (y)= 0$ and with $E_{A} (y^* y)= 1$.

Take $M$ to be any $II_1$-factor. We construct $M_1$ via the
following map $\Gamma$:
\[
  M_1 = \Gamma (M) = (M \otimes A)*_A M_{2n} (\mathbb {C}).
\]
The trace $\tau$ on $M$ can be extended to $M_1$.

Apply the above lemma.
\begin {align*}
  &P_1 = M \otimes A;\\
  &P_2 = M_{2n} (\mathbb {C});\\
  &L= P= M_1.
\end {align*}
We get:
\[
  {M_1}' \cap M_1= {M_1}' \cap A \subseteq {M_{2n} (\mathbb {C})}'
  \cap A = \mathbb {C}
\]
That is, $M_1$ is a nonhyperfinite $II_1$ factor.

The relative commutant
\[
  M' \cap M_1= (M' \cap M) \otimes A= A
\]
It also gives that $M_1' \cap M_1= {\mathbb C}$.

Viewing $\Gamma$ as a machine producing $II_1$ factors, we get an
ascending towers of $II_1$ factors:
\[
  M \subset \Gamma (M)= M_1 \subset \cdots M_i
  \subset \Gamma (M_i)= M_{i+1} \cdots.
\]
There is a unique trace $\tau$ associated to each $M_i$.
%
%
%

We calculate the relative commutant algebra $M' \cap M_k$ by
induction.
Let us assume
\[
  M' \cap M_k= \otimes^k A.
\]
By the above lemma,
\begin {align*}
  &M' \cap M_{k+1}= (M' \cap M_k) \otimes A\\
  &= (\otimes^k A) \otimes A= \otimes^{k+1} A.
\end {align*}
%
%
%

In the end, we would boldly suggest an analogy between binary shifts
and free products with amalgamation as in \cite{mB}.



\begin {acknowledgments}
I wish to express gratitude toward my advisor Vaughan F.R.Jones for
proposing the problem and also for useful comments and references. I
would also like to thank the referee for helping comments concerning
the organization of this paper.  A lot of helpful suggestions come
from Marie Choda, Zeph Landau, Geoffrey Price, and Erling
St{\o}rmer.
\end {acknowledgments}

\end {document}